\documentclass[twoside,a4paper,12pt,centertags,reqno]{amsart} 
\usepackage{amsmath,amssymb,verbatim,vmargin}
\usepackage{color}
\usepackage{tikz}
\usepackage{tikz-cd}
\usetikzlibrary{matrix,calc}
\usepackage{hyperref}
\hypersetup{colorlinks,bookmarks=true,linktocpage=true,citecolor=blue,linkcolor=magenta}

\usepackage{eulervm}   

\allowdisplaybreaks 
\usepackage{color}
\usepackage[all]{xy} 

\usepackage{mathtools}
\usepackage{MnSymbol} 

\theoremstyle{plain}
\newtheorem{thm}{Theorem}[section]

\theoremstyle{definition}

\newtheorem{rem}[thm]{Remark}

\usepackage{enumerate}
\usepackage{enumitem}

\newcommand{\bRn}{\mathbb{R}^n}

\newcommand{\pd}{\partial}

\newcommand{\bR}{{\mathbb R}}
\newcommand{\bS}{{\mathbb S}}

\newcommand{\fR}{{\mathbf R}}

\def\barint_#1{\mathchoice
            {\mathop{\vrule width 6pt
height 3 pt depth -2.5pt
                    \kern -9.5pt
\intop \kern -4pt}\nolimits_{#1}}%
            {\mathop{\vrule width 5pt height
3 pt depth -2.6pt
                    \kern -6.5pt
\intop \kern -4pt}\nolimits_{#1}}%
            {\mathop{\vrule width 5pt height
3 pt depth -2.6pt
                    \kern -6pt
\intop \kern -4pt}\nolimits_{#1}}%
            {\mathop{\vrule width 5pt height
3 pt depth -2.6pt
          \kern -6pt \intop \kern -4pt}\nolimits_{#1}}}
          
           \def\bariint_#1{\mathchoice
            {\mathop{\vrule width 15pt
height 3 pt depth -2.5pt
                    \kern -15.8pt
\intop \kern -8pt\intop \kern -4pt}\nolimits_{#1}}%
            {\mathop{\vrule width 9pt height
3 pt depth -2.6pt
                    \kern -10.5pt
\intop \kern -8pt\intop \kern -4pt}\nolimits_{#1}}%
            {\mathop{\vrule width 9pt height
3 pt depth -2.6pt
                    \kern -10pt
\intop \kern -8pt\intop \kern -4pt}\nolimits_{#1}}%
            {\mathop{\vrule width 9pt height
3 pt depth -2.6pt
          \kern -8pt \intop \kern -10pt\intop \kern -4pt}
      \nolimits_{  #1}}}

\def\barintlim_#1{\mathchoice
            {\mathop{\vrule width 6pt
height 3 pt depth -2.5pt
                    \kern -8.8pt
\intop \kern -4pt}\limits_{#1}}%
            {\mathop{\vrule width 5pt height
3 pt depth -2.6pt
                    \kern -6.5pt
\intop \kern -4pt}\limits_{#1}}%
            {\mathop{\vrule width 5pt height
3 pt depth -2.6pt
                    \kern -6pt
\intop \kern -4pt}\limits_{#1}}%
            {\mathop{\vrule width 5pt height
3 pt depth -2.6pt
          \kern -6pt \intop \kern -4pt}\limits_{#1}}}
          
           \def\bariintlim_#1{\mathchoice
            {\mathop{\vrule width 15pt
height 3 pt depth -2.5pt
                    \kern -15.8pt
\intop \kern -8pt\intop \kern -4pt}\limits_{#1}}%
            {\mathop{\vrule width 9pt height
3 pt depth -2.6pt
                    \kern -10.5pt
\intop \kern -8pt\intop \kern -4pt}\limits_{#1}}%
            {\mathop{\vrule width 9pt height
3 pt depth -2.6pt
                    \kern -10pt
\intop \kern -8pt\intop \kern -4pt}\limits_{#1}}%
            {\mathop{\vrule width 9pt height
3 pt depth -2.6pt
          \kern -8pt \intop \kern -10pt\intop \kern -4pt}
      \limits_{  #1}}}
          
\renewcommand{\iint}{\int \kern -8pt\int}       

\newcommand{\RE}{\text{{\rm Re}}\,}

\numberwithin{equation}{section}
\makeatletter
\@namedef{subjclassname@2020}{\textup{2020} Mathematics Subject Classification}
\makeatother

\title{Rellich Inequality via Radial Dissipativity}

\author{Yi C. Huang} 
\address{School of Mathematical Sciences, Nanjing Normal University, Nanjing 210023, People's Republic of China}
\email{Yi.Huang.Analysis@gmail.com}
\urladdr{https://orcid.org/0000-0002-1297-7674}

\date{\today} 

\thanks{Research of the author is supported by the National NSF grant of China (no. 11801274).
The author would like to thank Professor Tohru Ozawa (Waseda University) for kind support.}

\begin{document}

\begin{abstract}
We give a conceptually simple and essentially one-dimensional approach to Rellich inequality in Euclidean space $\bRn$.
In particular, we show that the radial part and the spherical part of the standard Laplacian form an angle in $[0,\pi/2]$ when $n\geq4$,
a property known in the works of Evans-Lewis, Machihara-Ozawa-Wadade, and Bez-Machihara-Ozawa.
Our proof here is direct and short.
\end{abstract}

\maketitle


\section{Introduction}

Let $\fR_n=\frac{n(n-4)}{4}$.
For proper $f$, we consider the following Rellich inequality:
\begin{equation} \label{e:Rellich}
\fR_n^2\left\|\frac{f}{|x|^2}\right\|_2^2\leq\|\Delta f\|_2^2,
\end{equation}
where $|\cdot|$ and $\|\cdot\|_2$ denote the Euclidean distance and the $L^2$-norm on $\bRn$.
Note that \eqref{e:Rellich} for $n=4$ is trivial,
and if one wishes to treat general functions in the Sobolev space $H^2(\bRn)$, 
it is necessary to restrict to $n\geq5$.
For the Laplacian, we have
$$\Delta=\Delta_r+\frac{1}{|x|^2}\Delta_{\bS^{n-1}},$$
where $\Delta_r$ denotes the radial Laplacian $\pd_r^2+\frac{n-1}{|x|}\pd_r$, with the radial derivative $\pd_r$ given by $\frac{x}{|x|}\cdot\nabla$,
and $\Delta_{\bS^{n-1}}$ denotes the Laplace-Beltrami operator on the sphere $\bS^{n-1}$.
Let
$$L_j=\pd_j-\frac{x_j}{|x|}\pd_r,\quad j=1, \cdots, n,$$
viewed as spherical-type derivatives.
Thus we have the following relation
$$\Delta_{\bS^{n-1}}=|x|^2\sum_{j=1}^nL_j^2.$$
Note that $-\sum_{j=1}^nL_j^2$ is a non-negative operator on $L^2(\bRn)$.

Recently in \cite{MacOzaWad17} and \cite{BezMacOza23}, 
Machihara-Ozawa-Wadade and Bez-Machihara-Ozawa took the perspective of proving the Rellich inequalities in the framework of equalities,
and found the following beautiful equalities: 
for $n\geq2$ and $f\in C_0^\infty(\bRn\backslash\{0\})$,
\begin{equation} \label{e:equal}
\|\Delta f\|_2^2=\left\|\Delta_r f\right\|_2^2+
\left\|\sum_{j=1}^nL_j^2 f\right\|_2^2+2\fR_n\sum_{j=1}^n\left\|\frac{L_jf}{|x|}\right\|_2^2+2\left\langle-\sum_{j=1}^nL_j^2f_*,f_*\right\rangle,
\end{equation}
where $f_*=\pd_rf+\frac{n-4}{2}\frac{f}{|x|}$.
Here $\langle\cdot,\cdot\rangle$ denotes the inner product for $L^2(\bRn)$.

The last two terms in \eqref{e:equal} are non-negative for $n\geq4$.
This property can be found in Evans-Lewis' continuation \cite{EvaLew05} of Laptev-Weidl type Hardy inequalities \cite{LapWei99}.
Our aim here is to expose a direct, simple, and short proof of this property that is also independent of the explicit representation in \eqref{e:equal} via $L_j$ and $f_*$.

For convenience, set 
$$\Delta_s=\frac{1}{|x|^2}\Delta_{\bS^{n-1}}=\sum_{j=1}^nL_j^2.$$

\begin{thm}[Evans-Lewis] \label{thm:radsph}
For $n\geq4$ and $f\in C_0^\infty(\bRn\backslash\{0\})$,
\begin{equation} \label{e:cross}
\RE\langle \Delta_r f, \Delta_sf\rangle\geq0,
\end{equation}
hence
\begin{equation} \label{e:Hilbert}
\|\Delta f\|^2_2\geq\|\Delta_r f\|^2_2+\|\Delta_s f\|^2_2.
\end{equation}
\end{thm}

\begin{rem}
As observed in \cite{BezMacOza23}, for $n=3$ one has
$$\|\Delta f\|^2_2\geq\|\Delta_r f\|^2_2\quad\text{and}\quad \|\Delta f\|^2_2\geq\left(\frac{5}{8}\right)^2\|\Delta_s f\|^2_2.$$
In this $n=3$ case it is interesting to have a refined angle description for \eqref{e:cross}.
\end{rem}

For more on Rellich inequalities, see the references in \cite{EvaLew05, MacOzaWad17, BezMacOza23}.
For related ideas on homogeneous groups, see Ruzhansky and Suragan \cite{RuzSur17, RuzSur19}.

\section{Proof of Theorem \ref{thm:radsph}}

Let $g$ be a function defined on $\bR_+$.
A simple integration-by-parts argument shows that $\frac{d}{dr}$ is a dissipative operator in $L^2(\bR_+, r^\alpha dr)$ for $\alpha\geq0$\footnotemark
\footnotetext{See for example Hayashi-Ozawa's revisit \cite[Section 3]{HayOza17} for Tosio Kato's insightful and far-reaching extension \cite{Kat71} 
of Landau-Kolmogorov inequalities to dissipative operators.}.
Hence, for $n\geq4$,
$$\RE\left\langle \frac{1}{r}\pd_rg,\frac{g}{r^2}\right\rangle_{L^2(\bR_+, r^{n-1}dr)}\leq0.$$
Also note that $\pd_r^2+\frac{(n-2)-1}{r}\pd_r$ is a dissipative operator in $L^2(\bR_+, r^{(n-2)-1}dr)$, i.e.,
$$\RE\left\langle \left(\pd_r^2+\frac{(n-2)-1}{r}\pd_r\right) g,\frac{g}{r^2}\right\rangle_{L^2(\bR_+, r^{n-1}dr)}\leq0,$$
seen as/via the dissipativity of the Laplacian on $\bR^{n-2}$ restricted to radial functions\footnotemark
\footnotetext{This change-of-dimension trick is useful in wave propagation and its application to Strichartz estimates in the non-radial regime, 
see for example Kato-Nakamura-Ozawa \cite[Lemma 2.1]{KatNakOza07}.}.
The above two inequalities via the expansion into spherical harmonics lead to \eqref{e:cross}.

\begin{rem}
Using polar coordinates and Fubini theorem, we see that
$$\|\Delta_r f\|^2_2\geq\fR_n^2\left\|\frac{f}{|x|^2}\right\|_2^2,\quad f\in C_0^\infty(\bRn\backslash\{0\}),$$
can be reduced to weighted Rellich inequalities on the half line, see \cite{HuaShi23}.
Combining this with \eqref{e:Hilbert} proposes an essentially one-dimensional approach to \eqref{e:Rellich}.
\end{rem}

\begin{rem}
All together we reduce the Rellich inequality on $\bRn$ to the dissipativity of $\frac{d}{dr}$ in $L^2(\bR_+, r^{n-4} dr)$ for $n\geq4$.
This is the essential part where $n\geq4$ plays a role.
\end{rem}

\bigskip

\section*{\textbf{Compliance with Ethical Standards}}

\bigskip

\textbf{Conflict of interest} The author has no known competing financial interests 
or personal relationships that could have appeared to influence this reported work.

\bigskip

\textbf{Availability of data and material} Not applicable.

\bigskip

\bibliographystyle{alpha}
 
\bibliography{Hua-RellichRadialDissipa}

\end{document}